\documentclass[A4paper,12pt]{article}
\usepackage{latexsym}
\usepackage{mathrsfs}
\usepackage{amssymb}
\usepackage{amscd}
\usepackage[dvips]{graphicx}                  

\newtheorem{theorem}{Theorem}

\newenvironment{proof}[1][Proof]{\textbf{#1.} }{\ \rule{0.5em}{0.5em}}
\date{}
\long\def\symbolfootnote[#1]#2{\begingroup%
	\def\thefootnote{$\;$}\footnote[#1]{$^*$#2}\endgroup}
\begin{document}
	
	\title{Weak $P$-points for regulars}
	\author{Joanna Jureczko}
\maketitle

\symbolfootnote[2]{Mathematics Subject Classification: 54D35, 54E45, 54C10.

	\hspace{0.2cm}
	}

\begin{abstract}
	The main result of this paper is to show that for any compact space $X$ and any family $\{f_\alpha \colon \alpha< 2^\kappa, f_\alpha \colon X \stackrel{onto}{\longrightarrow} \beta \kappa\}$ of open mappings  there exists a point  $x \in X$ such that for each $\alpha< 2^\kappa$ either $f_\alpha(x)$ belongs to $\kappa$ or $f_\alpha(x)$ is a weak $P_\kappa$-point in $\kappa^*$, whenever $\kappa$ is regular.
 	\end{abstract}

 \section{Introduction}
 
 The investigations of $P$-points and weak $P$-points have long history. The existence of these two notions are usually considered in $\omega^* = \beta\omega\setminus \omega$, where $\beta\omega$ means the \v Cech-Stone compactification. 
 Recall that Rudin in \cite{WR} showed that CH implies that there are $2^{2^\omega}$ $P$-points in $\omega^*$ and the same results holds under MA, but Shelah in \cite{SS}  showed that it is consistent with the axioms of set theory that there are no $P$-points in $\omega^*$. Furthermore, Kunen in \cite{KK} showed that there are $2^{2^\omega}$ points in $\omega^*$ which are weak $P$-points but not $P$-points. Moreover, Juh\'asz and Kunen in  \cite{JK} obtained two examples of such spaces. Additionally, in paper \cite{BK} Baker and Kunen generalized the results obtained in \cite{KK}. 
 
 On the other hand, in \cite{RF}  Frankiewicz obtains under some weaker version of MA that  for any compact space $X$ and any  family of open maps $$\{\pi_\xi \colon \xi < 2^\omega\}$$ with $\pi_\xi \colon X \stackrel{onto}{\longrightarrow} \beta\omega$  there exists a point  $x \in X$ such that for each $\xi < 2^\omega$  either $\pi_\xi(x)$ belongs to $\omega$ or $\pi_\xi(x)$ is not a limit point in $\omega^*$.
 This result was used for proving the following statements:
 \begin{enumerate}
 	\item In each extremally disconnected compact space there exists  a non-limit point of any countable discrete subset.
 	\item In $G(I)$ (Gleason space over the interval $[0,1]$) there exists a non-limit point of any discrete subset of cardinality $<2^\omega$.
 	\item Kunen's Theorem on non-limit points from \cite{KK1}. 
 \end{enumerate} 
Moreover, in \cite{JJ} the author obtained the above results under ZFC.
In this paper the main result of \cite{JJ} is generalized.  
For this purpose one need the method which allows us to continue the induction up to $2^\kappa$ where $\kappa> \omega$. With help it comes the results by Baker and Kunen published in \cite{BK} in which the authors  presented  very usefull method that can be recognized as a generalization of method  presented in \cite{KK}. 
It is worth emphasizing that both methods, (from \cite{KK} and \cite{BK}), provide usefull "technology"  for keeping the transfinite construction for an ultrafilter not finished before $\mathfrak{c}$ steps, (see \cite{KK}), and $2^\kappa, $  for $\kappa$ being infinite cardinal, (see \cite{BK}), but the second method has some limitations, among others $\kappa$ must be regular.

In the main results of this paper we will use a modification of arguments used in \cite[Theorem 6.1.]{BK} and some ideas used in \cite{RF}.

  The paper is divided into three sections. In Section 2 there are gathered main definitions and facts used in the further parts of the paper. Section 3 contains main results. Section 4 contains remarks and open questions.
 
 All notations used in the paper are standard for the area and are assumed as well known. For definitions and facts not cited here the reader is referred to~\cite{RE, FZ, TJ}.
 
 \section{Definitions and previous results}
 \textbf{2.1.} Let $\kappa$ be an infinite cardinal. Let $Y$ be a topological space and $p\in Y$. We say that
 \begin{itemize}
 	\item $p$ is a \textit{$P_{\kappa}$-point} in $Y$ iff the intersection of any family of fewer than $\kappa$ neighbourhoods of $p$ is also a neighbourhood of $p$,
 	\item $p$ is a \textit{weak $P_\kappa$-point} in $Y$ iff $p$ is not a limit point of any subset of $Y \setminus\{p\}$ of size less than $\kappa$.
 \end{itemize}
Obviously, a $P_{\omega_1}$-point  and a weak $P_{\omega_1}$-point are called a \textit{$P$-point} and a \textit{weak $P$-point}, respectively.
\\
\\
\textbf{2.2.} Let $\kappa$ be an infinite cardinal. Then $\beta\kappa$ means the \v Cech-Stone compactification, where $\kappa$ has the discrete topology. Hence, $\beta\kappa$ is the space of ultrafilters on $\kappa$ and $\kappa^* = \beta \kappa \setminus \kappa$ is the space of nonprincipal ultrafilters on $\kappa$. Let $\mathcal{U}(\kappa)$ denote the space of uniform ultrafilters, i.e. the space of ultrafilters of size $\kappa$. Obviously, $\mathcal{U}(\kappa) \subseteq \kappa^* \subseteq \beta \kappa$. 
\\\\
\textbf{2.3.} 	Let us accept the following notation: 
\begin{itemize}
	\item $\mathcal{FR}(\kappa) = \{A \subset \kappa \colon |\kappa\setminus A|< \kappa\},$
	\item $[A, B, C,...]$ means the filter generated by $A, B, C, ...$,
	\item $A\subseteq_* B$ iff $|A\setminus B|< \kappa$  for any $A, B \subseteq \kappa$.
\end{itemize} 
\textbf{2.4.} A function $\hat\varphi \colon [\kappa^+]^{<\omega} \to [\kappa]^{< \omega}$ is \textit{$\kappa$-shrinking} iff 
\begin{itemize}
	\item $p\subseteq q$ implies $\hat{\varphi}(p) \subseteq \hat{\varphi}(q)$, for any $p, q \in [\kappa^+]^{<\omega}$,
	\item $\hat{\varphi}(0) = 0$.
\end{itemize}

A \textit{step family} (over $\kappa$, with respect to $\hat{\varphi}$) is a family of subsets of $\kappa$, $$\{E_t \colon t \in [\kappa]^{<\omega}\} \cup \{A_\alpha \colon \alpha < \kappa^+\}$$ satisfying the following conditions:
\begin{itemize}
	\item $E_s\cap E_t = \emptyset$ for all $s, t \in [\kappa]^{<\omega}$ with $s \not =t$,
	\item $|\bigcap_{\alpha \in p} A_\alpha \cap  \bigcup_{t\not \supseteq \hat{\varphi}(p)}E_t| < \kappa$ for each $p \in [\kappa^+]^{<\omega}$,
	\item if $\hat{\varphi}(p) \subseteq t$, then $|\bigcap_{\alpha\in p} A_\alpha \cap E_t|=\kappa$ for each $p \in [\kappa^+]^{<\omega}$ and $t \in [\kappa]^{<\omega}$. 
\end{itemize}

Let $I$ be an index set and $\mathcal{F}$ be a filter on $\kappa$. The family $$\{E_t^i \colon t \in [\kappa]^{<\omega}, i \in I\} \cup \{A_\alpha^i \colon \alpha < \kappa^+, i \in I\}$$ is an \textit{independent matrix of $|I|$ step-families} (over $\kappa$) with respect to $\mathcal{F}, \hat{\varphi}$ iff
\begin{itemize}
	\item for each fixed $i \in I$, $\{E_t^i \colon t \in [\kappa]^{<\omega}\} \cup \{A_\alpha^i \colon \alpha < \kappa^+\}$ is a step-family,
	\item if $n \in \omega, p_0, p_1, ..., p_{n-1} \in [\kappa^+]^{<\omega}, t_0, t_1, ..., t_{n-1} \in [\kappa]^{<\omega}$,  $i_0, i_1, ..., i_{n-1} \in I$ with $i_k\not = i_m, k\not = m$ and $\hat{\varphi}(p_k) \subseteq t_k$, then 
	$$\bigcap_{k=1}^{n-1}(\bigcap_{\alpha\in p_k}A^{i_k}_{\alpha} \cap E^{i_k}_{t_k}) \in \mathcal{F}^+,$$
	where $\mathcal{F}^+ = \{D\subseteq \kappa \colon \kappa\setminus D \not \in \mathcal{F}\}$.  
\end{itemize}
\noindent
\textbf{Fact 1 (\cite{BK})} If $\kappa$ is a regular cardinal and $\hat{\varphi}$ is a $\kappa$-shrinking function, then there exists and independent matrix of $2^\kappa$ step-families over $\kappa$ with respect to $\mathcal{FR}(\kappa)$, $\hat{\varphi}$. 
\\\\
\textbf{2.5.}
Let $\hat{\varphi}$ be a $\kappa$-shrinking function and let $x \in Y$. A point $x$ is called a \textit{$\kappa$-shrinking point} iff whenever $\{U_t \colon t \in [\kappa]^{<\omega}\}$ are neighourhoods of $x$ then there exist neighbourhoods $\{V_\beta \colon \beta < \kappa^+\}$ of $x$ such that $\bigcap_{\alpha \in p} V_\alpha \subseteq U_{\hat{\varphi}(p)}$ for each non-empty $p \in [\kappa^+]^{<\omega}$.
\\
\\
\textbf{Fact 2. (\cite{BK}).} Let $\kappa$ be regular and $Y$ be a space. There exists a $\kappa$-shrinking function such that each $\kappa$-shrinking point $x$ in $Y$ is a weak $P_\kappa$-point in $Y$.
\\\\
\textbf{Fact 3. (\cite{BK}).} Let $\kappa$ be a regular cardinal. There exists a $\kappa$-shrinking function such that whenever $x \in \mathcal{U}(\kappa)$ is a $\kappa$-shrinking point in $\mathcal{U}(\kappa)$, then $x$ is a $\kappa$-shrinking point in $\beta \kappa$. 
\\
\\
Notice that a $\kappa$-shrinking function and a $\kappa$-shrinking point were introduced in \cite{BK} as a hat function and a hat point.

\section{Main results}

This section contains four results. The first two concern the considerations of the existing of $\kappa$-shrinking points for one function, while the last two concern the considerations of existing $\kappa$-shrinking points which are common for $2^\kappa$ many functions. 
Since proofs of these results are essentially different one decided to present all of them with all details.
\\

 Let $X$ be a space, (since now, we assume that $X$ is a $T_2$-space), and let $$\{f_\alpha \colon f_\alpha \colon X \stackrel{onto}{\longrightarrow} \beta \kappa, \alpha < 2^\kappa\}$$ be a family of open mappings.

Let $\{S_\beta \colon \beta < 2^\kappa\}$ be a family of non-empty and closed subsets of $X$. For each $\gamma < 2^\kappa$ let 
$$C_\gamma = \bigcap\{S_\beta \colon \beta < \gamma\}$$
and
$$B_\gamma = \left\{\begin{array}{rcl}
C_\gamma \cap f^{-1}_\nu(\kappa^*) & \textrm{where } \nu = \min \{\beta < 2^\kappa \colon C_\gamma \cap f^{-1}_\beta(\kappa^*) \not = \emptyset\}\\
C_\gamma & \textrm{otherwise}
\end{array} \right.$$

In Theorem 1 we will consider one function so we have to slightly modify the definition of $B_\gamma$, i.e. $B_\gamma = C_\gamma \cap f^{-1}(\kappa^*)$ whenever $C_\gamma \cap f^{-1}(\kappa^*) \not = \emptyset$.

\begin{theorem}
		Let $\kappa$ be a regular cardinal and let $\hat{\varphi}$  be a $\kappa$-shrinking function.
	Let $X$ be a compact space and let $ f \colon X \stackrel{onto}{\longrightarrow} \beta \kappa$ be an open mapping. Then  there exists a family $\{U_{\beta} \colon \beta < 2^\kappa\}$ of non-empty and closed subsets of $X$ such that for each $\gamma < 2^\kappa$
\begin{itemize}
	\item [(i)] $U_{\beta} \subseteq U_{\beta'}$ for any $\beta' \leqslant \beta < \gamma$;
	\item [(ii)] if $f^{-1}(\kappa^*) \cap B_\gamma = \emptyset$, then
	$U_{\gamma} = f^{-1}(\{\xi\}) \cap B_\gamma$ for some $\xi \in \kappa$;
	\item [(iii)] if $f^{-1}(\kappa^*) \cap B_\gamma  \not = \emptyset$, then there exists an ultrafilter $\bigcup_{\beta \in 2^\kappa}\mathcal{F}_\beta$ such that 
	\begin{itemize}
		\item [(a)] $\mathcal{F}_{\beta} \subseteq \mathcal{F}_{\beta'}$, whenever $\beta \leqslant \beta'< \gamma$
		\item [(b)] for each decreasing sequence $\{(V^r_\beta \subseteq \omega \colon r < [\kappa]^{<\omega}), \beta \in 2^\kappa\}$ if $V^r_\beta \in \mathcal{F}_\beta$ for $r < [\kappa]^{<\omega}$ then there exists $W^\eta_\beta \in \mathcal{F}_{\beta+1}$ for each $\eta \in \kappa^+$ such that for each $p \in [\kappa^+]^{<\omega}$ and each $\eta_1< \eta_2< ... < \eta_k < 2^\kappa$
		$$|\bigcap_{\nu \in p}W^{\nu}_{\beta} \setminus V^{\hat{\varphi}(p)}_{\beta}|< \kappa,$$
		\item [(c)] $f^{-1}(U) \cap B_\gamma \not =\emptyset$, where $U$ is a clopen set belonging to the standard base of $\kappa^*$ such that $\bigcup_{\beta< \gamma} \mathcal{F}_\beta \in U$.
	\end{itemize}	
\end{itemize}
\end{theorem}

\begin{proof}	
By Fact 1, fix a matrix $$\{E_t^i \colon t \in [\kappa]^{<\omega}, i \in 2^\kappa\} \cup \{A_\eta^i \colon \eta < \kappa^+, i \in 2^\kappa\}$$ of $2^\kappa$ step-families over $\kappa$ independent with respect to  $\mathcal{FR}(\kappa), \hat{\varphi}$. 
	Assume that for each $i \in 2^\kappa$ we have 
	\begin{itemize}
		\item [(a)] $E^i_s\cap E^i_t = \emptyset$ for all $s, t \in [\kappa]^{<\omega}$ with $s \not =t$,
		\item [(b)] $|\bigcap_{\eta \in p} A^i_\eta \cap  \bigcup_{t\not \supseteq \hat{\varphi}(p)}E^i_t| < \kappa$ for each $p \in [\kappa^+]^{<\omega}$,
		\item [(c)] if $\hat{\varphi}(p) \subseteq t$, then $|\bigcap_{\eta\in p} A^i_\eta \cap E^i_t|=\kappa$ for each $p \in [\kappa^+]^{<\omega}$ and $t \in [\kappa]^{<\omega}$ 
		\item [(d)] $\bigcup\{E^i_t \colon t \in [\kappa]^{< \omega}\} = \kappa$
		\item [(e)] $|\bigcap_{\eta \in p} A^i_\eta \setminus  \bigcup_{t \supseteq \hat{\varphi}(p)}E^i_t| < \kappa$ for each $p \in [\kappa^+]^{<\omega}$.
	\end{itemize} 
Note that $(b)$ and $(d)$ implies $(e)$.
Moreover, the condition $(b)$ is still preserved after expanding $\{E_t^i \colon t \in [\kappa]^{<\omega}\}$ to a partition of $\kappa$. 
\\Indeed. If there are $p_0 \in [\kappa^+]^{<\omega}$ and $i_0 \in 2^\kappa$ such that
$$|\bigcap_{\eta \in p_0} A_\eta^{i_0} \cap  \bigcup_{t\not \supseteq \hat{\varphi}(p_0)}E_t^{i_0}| = \kappa,$$
then $|\bigcup_{t\not \supseteq \hat{\varphi}(p_0)}E_t^{i_0}| = \kappa.$
Then, by $(d)$ and $(a)$, there would exist  $t_0 \supseteq \hat{\varphi}(p_0)$ such that  $|E_{t_0}^{i_0}| < \kappa$. Hence
$$|\bigcap_{\eta\in p_0} A^{i_0}_\eta \cap E^{i_0}_{t_0}|<\kappa.$$
which contradicts $(c)$.

Let $\{Z_\beta \colon \beta < 2^\kappa, \alpha \equiv 0\ (mod\  2)\}$ be the family of all  subsets of $\kappa$ and let $$\{(V^\beta_r \colon r < [\kappa]^{<\omega}) \colon \beta < 2^\kappa, \beta \equiv 1\ (mod\  2)\}$$ be the family of  monotone sequences, i.e. $V^\beta_r \supseteq V^\beta_s$ whenever $r \subseteq s$. 
 (We do not require the $(V^\beta_r\colon r < [\kappa]^{<\omega})$ are to be distinct).

We construct $\mathcal{F}_\beta$, $I_{\beta}$ and $U_{\beta}$ by induction in $2^\kappa$ steps in the following way, ($B_\beta$ are constructed as shown before the theorem but with using sets $U_\beta$) 
	\begin{itemize}
	\item[(1)] $\mathcal{F}_0 = \mathcal{FR}(\kappa)$, $I_0= 2^\kappa$ and $U_{0} = X$,
	\item[(2)] $\mathcal{F}_\beta$ is a filter on $\kappa$, $I_\beta \subset 2^\kappa$, $U_{\beta}$ is a non-empty closed subset of $X$ and $$\{E_t^i \colon t \in [\kappa]^{<\omega}, i \in I_\beta\} \cup \{A_\mu^i \colon \mu < \kappa^+, i \in I_\beta\}$$ of remaining step-families is independent with respect to $\mathcal{F}_\beta, \hat{\varphi}$,
	\item[(3)] if $\gamma < \beta$, then $\mathcal{F}_\gamma \subseteq \mathcal{F}_\beta$, $I_\gamma\supseteq I_\beta$, $U_{\gamma} \supseteq U_{\beta}$,
	\item[(4)] if $\beta$ is limit, then $\mathcal{F}_\beta = \bigcup_{\gamma< \beta} \mathcal{F}_\gamma$, $I_\beta = \bigcap_{\gamma<\beta} I_\gamma$, $U_{\beta} = \bigcap_{\gamma<\beta} U_{\gamma}$,
	\item[(5)] $I_\beta \setminus I_{\beta+1}$ is finite for each $\beta \in 2^\kappa$,
	\item[(6)] either $Z_\beta \in \mathcal{F}_\beta$ or $\kappa\setminus Z_\beta \in \mathcal{F}_\beta$, 
	\item[(7)] if $f(B_{\beta}) \cap \kappa^* = \emptyset$ then $U_{\beta} = B_{\beta} \cap f^{-1}(\{\xi\})$ for some $\xi \in \kappa$,
	\item[(8)] if $f(B_{\beta}) \cap \kappa^* \not = \emptyset$ then 
	\begin{itemize}
		\item  $f(B_{\beta}) \cap U \not =\emptyset$, for some clopen set $U$ belonging to the standard base of $\kappa^*$ and $\bigcup_{\gamma < \beta} \mathcal{F}_\gamma \in U$,
		\item  $U_{\beta} = f^{-1}(U) \subseteq B_{\beta}$
		\item  for a decreasing sequence $\{V_\beta^r \colon r \in [\kappa]^{< \omega}\}$ if each $V^r_\beta \in \mathcal{F}_\beta$ then there are $W^\beta_\eta \in \mathcal{F}_{\beta+1}$, for each $\eta \in \kappa^+$, such that for each $p \in [\kappa^+]^{<\omega}$, $$|\bigcap_{\nu \in p}W_\beta^\nu\setminus V_\beta^{\hat{\varphi}(p)}|< \kappa.$$
	\end{itemize}
\end{itemize} 	

The first and the limit step are done. Assume that we have constructed $\mathcal{F}_\beta$, $I_\beta$ and $U_{\beta}$. We will show the successor step.  

Assume that  $f(B_{\beta}) \cap \kappa^* \not = \emptyset$. Otherwise, by $(7)$ we put $U_{\beta+1} = B_{\beta+1} \cap f^{-1}(\{\xi\})$ for some $\xi \in \kappa$. Then $(ii)$ holds.

Consider two cases.

Case 1. $\beta \equiv 0 \ (mod\  2)$. Let  $\mathcal{R} = [\mathcal{F}_\beta, \{Z_\beta\}]$ be a filter.
If $\mathcal{R}$ is proper and 
$$\{E_t^i \colon t \in [\kappa]^{<\omega}, i \in I_\beta\} \cup \{A_\mu^i \colon \mu < \kappa^+, i \in I_\beta\}$$ is independent with respect to $\mathcal{R}, \hat{\varphi}$, then  put 
then we put
$$\mathcal{F}_{\beta+1} = \mathcal{F}_\beta, \ \ I_{\beta+1} = I_\beta, \ \ U_{\beta+1} = U_{\beta}.$$

Otherwise, fix $n \in \omega$, distinct $i_k \in I_\beta$ and $\hat{\varphi}(p_k)\subseteq t_k,$ for $ k < n$, such that
$$\kappa \setminus [Z_\beta \cap \bigcap_{k=0}^{n-1}(\bigcap_{\mu\in p_k}A^{i_k}_{\mu} \cap E^{i_k}_{t_k})]\in \mathcal{F}_\beta.$$
Then, put $$\mathcal{F}_{\beta+1} = [\mathcal{F}_{\beta}, \{A^{i_k}_{p_k} \colon 0 \leqslant k \leqslant n-1\}, \{E^{i_k}_{t_k} \colon 0 \leqslant k \leqslant n-1\}]$$ $$I_{\beta +1} = I_{\beta}\setminus \{i_k \colon 0 \leqslant k \leqslant n-1\}$$
and
$$U_{\beta +1} = f^{-1}(U) \subseteq B_{\beta+1},$$
where $U$ is a clopen set belonging to the standard basis of $\kappa^*$ and $\bigcup_{\gamma<\beta+1}\mathcal{F}_\gamma \in U$.
Then $\kappa\setminus Z_\beta \in \mathcal{F}_{\beta+1}.$

To show that $(2)$ is fulfilled for $\mathcal{F}_{\beta+1}$ and $I_{\beta+1}$  it is enough to observe that each element of $\mathcal{F}_{\beta+1}$ is of the form 
$$A \cap (\bigcap_{\mu \in p_k}A^{i_k}_\mu \cap E^{i_k}_{t_k})$$
for some $A \in \mathcal{F}_{\beta} $ and $ 0 \leqslant k \leqslant n-1$. Then,
$$\bigcap_{\mu \in p_k}A^{i_k}_\mu \cap E^{i_k}_{t_k} \supseteq\bigcap_{k=0}^{n-1}(\bigcap_{\mu\in p_k}A^{i_k}_{\mu} \cap E^{i_k}_{t_k}) $$
and $$\bigcap_{k=0}^{n-1}(\bigcap_{\mu\in p_k}A^{i_k}_{\mu} \cap E^{i_k}_{t_k}) \in \mathcal{F}_{\beta}.$$
Thus, $(2)$ for $\mathcal{F}_{\beta+1}$ and $I_{\beta+1}$ is fulfilled which follows from $(2)$ for $\mathcal{F}_{\beta}$ and $I_\beta$. 

Case 2. $\beta \equiv 1 \ (mod\  2)$.

If $V^r_\beta \not \in \mathcal{F}_\beta$ for some $r \in [\kappa]^{<\omega}$ then we put
$$\mathcal{F}_{\beta+1} = \mathcal{F}_\beta, \ \ I_{\beta+1} = I_\beta, \ \ U_{\beta+1} = U_{\beta}.$$

Assume that each $V^r_\beta \in \mathcal{F}_\beta$. Choose $i \in I_\beta$. Then for each $\nu < \kappa^+$ consider
$$W_\beta^\nu = A^i_\nu \cap \bigcup\{V^t_\beta \cap E^i_t \colon t \in [\kappa]^{< \omega}\}.$$
Then for each $p \in [\kappa^+]^{<\omega}$ we have
$$\bigcap_{\nu\in p} W_\beta^\nu = \bigcap_{\nu\in p} A^i_\nu \cap \bigcup\{V^t_\beta \cap E^i_t \colon t \in [\kappa]^{< \omega}\} \subseteq_*$$
then by $(e)$
$$\bigcup\{V^t_\beta \cap E^i_t \colon t \in [\kappa]^{< \omega}\} \cap \{E^i_t \colon t \supseteq \hat{\varphi}(p)\} =$$
$$\bigcup \{V^t_\beta \cap E^i_t \colon t \supseteq \hat{\varphi}(p)\} \subseteq $$
$$\{V^t_\beta \colon t \supseteq \hat{\varphi}(p)\} = \bigcap_{\nu \in \hat{\varphi}(p)}V^\nu_\beta.$$
The last equality is obtained by the monotonicity of $\{V^r_\beta \colon r \in [\kappa]^{<\omega}\}.$

Now, put 
$$\mathcal{F}_{\beta+1} = [\mathcal{F}_\beta, \{W^\nu_\beta \colon \nu< \kappa^+\}]$$
$$I_{\beta+1} = I_\beta \setminus \{i\}$$
and
$$U_{\beta +1} = f^{-1}(U) \subseteq B_{\beta+1},$$
where $U$ is a clopen set belonging to the standard basis of $\kappa^*$ and $\bigcup_{\gamma < \beta+1}\mathcal{F}_\gamma \in U$.

To show that $(2)$ holds for $\mathcal{F}_{\beta+1}$ and $I_{\beta+1}$ it is enough to observe that each element of $\mathcal{F}_{\beta+1}$  is of the form 
$$A \cap \bigcap_{\nu \in p}W_\beta^\nu$$
for some $A\in \mathcal{F}_\beta$ and $p \in [\kappa^+]^{<\omega}$.
But 
$$V^{\hat{\varphi}(p)}_{\beta} \cap \bigcap_{\nu \in p} A^i_\nu\cap E^i_{\hat{\varphi}(p)} \subseteq \bigcap_{\nu \in p}W^\nu_\beta$$
and since each $V^r_\beta \in \mathcal{F}_\beta$ we have that $V^{\hat{\varphi}(p)}_{\beta} \in \mathcal{F}_\beta$.
Thus, $(2)$ for $\mathcal{F}_{\beta+1}$ and $I_{\beta+1}$ is fulfilled which follows from $(2)$ for $\mathcal{F}_\beta$ and $I_\beta$. 
\end{proof}
\\

\begin{theorem}
		Let $\kappa$ be a regular cardinal and let $\hat{\varphi}$  is a $\kappa$-shrinking function.
	Let $X$ be a compact space and let $f \colon X \stackrel{onto}{\longrightarrow} \beta \kappa$ be an open mapping. Then there exists a family of sets $\{U_\beta \colon \beta < 2^\kappa\}$  fulfilling $(i)-(iii)$ of Theorem 1 such that for each  $x \in \bigcap\{U_\beta \colon \beta < 2^\kappa\}$  either $f(x)$ belongs to $\kappa$ or $f(x)$ is a $\kappa$-shrinking point in $\kappa^*$.
\end{theorem}

\begin{proof}
	By Theorem 1, there exists a family $\{U_{\beta} \colon \beta \in 2^\kappa\}$ of non-empty and closed subsets of $X$ and an ultrafilter $\bigcup_{\beta \in 2^\kappa}\mathcal{F}_\beta$ of properties $(i)-(iii)$. 
	Then for each $x \in \bigcap\{U_{\beta} \colon \beta \in 2^\kappa\}$ $f(x)$ belongs to $\kappa$ whenever
	$f_{\alpha}(B_{\beta}) \cap \kappa^* = \emptyset$ or $f(x)$    is a $\kappa$-shrinking point in $\kappa^*$
	(of the form $\bigcup_{\beta \in 2^\kappa}\mathcal{F}_{\beta}$) whenever 
	$ f(B_{\beta}) \cap \kappa^* \not = \emptyset$.
\end{proof}
\\

Let $<_*$ denote the canonical well-ordering on $2^\kappa \times 2^\kappa$. Let $j(\alpha, \beta)\in 2^\kappa \times 2^\kappa$ be an immediate successor of $(\alpha, \beta)$, i.e. $(\alpha, \beta)<_*j(\alpha, \beta)$ and there is no $(\gamma, \delta)$ such that $(\alpha, \beta)<_*(\gamma, \delta)<_*j(\alpha, \beta)$.

Let $X$ be a space and let $\{f_\alpha \colon f_\alpha \colon X \stackrel{onto}{\longrightarrow} \beta \kappa, \alpha < 2^\kappa\}$ be a family of open mappings. 
Let $\{S_{\alpha, \beta} \colon (\alpha, \beta) < 2^\kappa\times 2^\kappa\}$ be a family of non-empty and closed subsets of $X$. For each $(\gamma, \delta) \in 2^\kappa \times 2^\kappa$ let 
$$C_{\gamma, \delta} = \bigcap\{S_{\alpha,\beta} \colon (\alpha, \beta) <_* (\gamma, \delta)\}$$ and
$$B_{\gamma, \delta} = \left\{\begin{array}{rcl}
C_{\gamma, \delta} \cap f^{-1}_\mu(\kappa^*) & \textrm{where } \mu = \min \{\beta < 2^\kappa \colon C_{\gamma, \delta} \cap f^{-1}_\beta(\kappa^*) \not = \emptyset\}\\
C_{\gamma, \delta} & \textrm{otherwise}
\end{array} \right.$$

\begin{theorem}
	Let $\kappa$ be a regular cardinal and let $\hat{\varphi}$  be a $\kappa$-shrinking function.
	Let $X$ be a compact space and let $\{f_\alpha \colon f_\alpha \colon X \stackrel{onto}{\longrightarrow} \beta \kappa, \alpha < 2^\kappa\}$ be a family of open mappings. Then  there exists a family $\{U_{\alpha,\beta} \colon (\alpha, \beta) < 2^\kappa\times 2^\kappa\}$ of non-empty and closed subsets of $X$ such that for each $(\gamma, \delta) < 2^\kappa\times 2^\kappa$
	\begin{itemize}
		\item [(i)] $U_{\alpha, \beta} \subseteq U_{\alpha, \beta'}$ for any $(\alpha',\beta') \leqslant_* (\alpha,\beta) <_* (\gamma, \delta)$;
		\item [(ii)] if $f^{-1}_\gamma(\kappa^*) \cap B_{\gamma, \delta} = \emptyset$, then
		$U_{\gamma, \delta} = f^{-1}_\gamma(\{\xi\}) \cap B_{\gamma, \delta}$ for some $\xi \in \kappa$;
		\item [(iii)] if $f^{-1}_\gamma(\kappa^*) \cap B_{\gamma, \delta}  \not = \emptyset$, then there exists an ultrafilter $\bigcup_{(\alpha, \beta) \in 2^\kappa\times 2^\kappa}\mathcal{F}_{\alpha, \beta}$ such that 
		\begin{itemize}
			\item [(a)] $\mathcal{F}_{\alpha, \beta} \subseteq \mathcal{F}_{\alpha', \beta'}$, whenever $(\alpha,\beta) \leqslant_* (\alpha',\beta') <_* (\gamma, \delta)$;
			\item [(b)] for each decreasing sequence $\{(V^r_{\alpha, \beta} \subseteq \omega \colon r < [\kappa]^{<\omega}), (\alpha,\beta) \in 2^\kappa\times 2^\kappa\}$ if $V^r_{\alpha, \beta} \in \mathcal{F}_{\alpha, \beta}$ for $r < [\kappa]^{<\omega}$ then there exists $W^\eta_{\alpha, \beta} \in \mathcal{F}_{j(\alpha, \beta)}$ for each $\eta \in \kappa^+$ such that for each $p \in [\kappa^+]^{<\omega}$ and each $\eta_1< \eta_2< ... < \eta_k < 2^\kappa$
			$$|\bigcap_{\nu \in p}W^{\nu}_{\alpha, \beta} \setminus V^{\hat{\varphi}(p)}_{\alpha, \beta}|< \kappa,$$
			\item [(c)] $f_\gamma^{-1}(U) \cap B_{\gamma, \delta} \not =\emptyset$, where $U$ is a clopen set belonging to the standard base of $\kappa^*$ such that $\bigcup_{(\alpha, \beta)<_* (\gamma, \delta)} \mathcal{F}_{\alpha, \beta} \in U$.
		\end{itemize}	
	\end{itemize}
\end{theorem}

\begin{proof}	
	By Fact 1, fix a matrix $$\{E_t^i \colon t \in [\kappa]^{<\omega}, i \in 2^\kappa\} \cup \{A_\eta^i \colon \eta < \kappa^+, i \in 2^\kappa\}$$ of $2^\kappa$ step-families over $\kappa$ independent with respect to $\mathcal{FR}(\kappa), \hat{\varphi}$. 
	Assume that for each $i \in 2^\kappa$ we have 
	\begin{itemize}
		\item [(a)] $E^i_s\cap E^i_t = \emptyset$ for all $s, t \in [\kappa]^{<\omega}$ with $s \not =t$,
		\item [(b)] $|\bigcap_{\eta \in p} A^i_\eta \cap  \bigcup_{t\not \supseteq \hat{\varphi}(p)}E^i_t| < \kappa$ for each $p \in [\kappa^+]^{<\omega}$,
		\item [(c)] if $\hat{\varphi}(p) \subseteq t$, then $|\bigcap_{\eta\in p} A^i_\eta \cap E^i_t|=\kappa$ for each $p \in [\kappa^+]^{<\omega}$ and $t \in [\kappa]^{<\omega}$ 
		\item [(d)] $\bigcup\{E^i_t \colon t \in [\kappa]^{< \omega}\} = \kappa$
		\item [(e)] $|\bigcap_{\eta \in p} A^i_\eta \setminus  \bigcup_{t \supseteq \hat{\varphi}(p)}E^i_t| < \kappa$ for each $p \in [\kappa^+]^{<\omega}$.
	\end{itemize} 
	Note that $(b)$ and $(d)$ implies $(e)$.
	Moreover, the condition $(b)$ is still preserved after expanding $\{E_t^i \colon t \in [\kappa]^{<\omega}\}$ to a partition of $\kappa$. 
	\\Indeed. If there are $p_0 \in [\kappa^+]^{<\omega}$ and $i_0 \in 2^\kappa$ such that
	$$|\bigcap_{\eta \in p_0} A_\eta^{i_0} \cap  \bigcup_{t\not \supseteq \hat{\varphi}(p_0)}E_t^{i_0}| = \kappa,$$
	then $|\bigcup_{t\not \supseteq \hat{\varphi}(p_0)}E_t^{i_0}| = \kappa.$
	Then, by $(d)$ and $(a)$, there would exist  $t_0 \supseteq \hat{\varphi}(p_0)$ such that  $|E_{t_0}^{i_0}| < \kappa$. Hence
	$$|\bigcap_{\eta\in p_0} A^{i_0}_\eta \cap E^{i_0}_{t_0}|<\kappa.$$
	which contradicts $(c)$.

	Let $\{Z_{\alpha, \beta} \colon (\alpha, \beta) < 2^\kappa\times 2^\kappa, \alpha + \beta \equiv 0\ (mod\  2)\}$ be the family of all  subsets of $\kappa$ and let $$\{(V_{\alpha, \beta}^r \colon r < [\kappa]^{<\omega}) \colon (\alpha, \beta) < 2^\kappa\times 2^\kappa, \alpha+\beta \equiv 1\ (mod\  2)\}$$ be the family of  monotone sequences, i.e. $V_{\alpha, \beta}^r \supseteq V_{\alpha, \beta}^s$ whenever $r \subseteq s$. 
	(We do not require the $(V_{\alpha, \beta}^r\colon  r < [\kappa]^{<\omega}))$ are to be distinct).
	
	We construct $\mathcal{F}_{\alpha, \beta}$, $I_{\alpha, \beta}$ and $U_{\alpha, \beta}$ by induction in $2^\kappa$ steps in the following way, ($B_{\alpha, \beta}$ will be defined as is given before the theorem),
	\begin{itemize}
		\item[(1)] $\mathcal{F}_{0,0} = \mathcal{FR}(\kappa)$, $I_{0,0}= 2^\kappa$ and $U_{0, 0} = X$,
		\item[(2)] $\mathcal{F}_{\alpha, \beta}$ is a filter on $\kappa$, $I_{\alpha, \beta} \subset 2^\kappa$, $U_{\alpha, \beta}$ is a non-empty closed subset of $X$ and $$\{E_t^i \colon t \in [\kappa]^{<\omega}, i \in I_{\alpha, \beta}\} \cup \{A_\mu^i \colon \mu < \kappa^+, i \in I_{\alpha, \beta}\}$$ of remaining step-families is independent with respect to $\mathcal{F}_{\alpha, \beta}, \hat{\varphi}$,
		\item[(3)] if $(\gamma, \delta) <_* (\alpha, \beta)$, then $\mathcal{F}_{\gamma, \delta} \subseteq \mathcal{F}_{\alpha, \beta}$, $I_{\gamma, \delta}\supseteq I_{\alpha, \beta}$, $U_{\gamma, \delta} \supseteq U_{\alpha, \beta}$,
		\item[(4)] if $(\alpha, \beta)$ is limit, then $\mathcal{F}_{\alpha, \beta} = \bigcup_{(\gamma, \delta)<_* (\alpha, \beta)} \mathcal{F}_{\gamma, \delta}$, $I_{\alpha, \beta} = \bigcap_{(\gamma, \delta)<_*(\alpha, \beta)} I_{\gamma, \delta}$, $U_{\alpha, \beta} = \bigcap_{(\gamma, \delta)<_*(\alpha, \beta)} U_{\gamma, \delta}$,
		\item[(5)] $I_{\alpha, \beta} \setminus I_{j(\alpha, \beta)}$ is finite for each $(\alpha, \beta) \in 2^\kappa\times 2^\kappa$,
		\item[(6)] either $Z_{\alpha, \beta} \in \mathcal{F}_{\alpha, \beta}$ or $\kappa\setminus Z_{\alpha, \beta} \in \mathcal{F}_{\alpha, \beta}$, 
		\item[(7)] if $f_\alpha(B_{\alpha, \beta}) \cap \kappa^* = \emptyset$ then $U_{\alpha, \beta} = B_{\alpha, \beta} \cap f_\alpha^{-1}(\{\xi\})$ for some $\xi \in \kappa$,
		\item[(8)] if $f_\alpha(B_{\alpha, \beta}) \cap \kappa^* \not = \emptyset$ then 
		\begin{itemize}
			\item  $f_\alpha(B_{\alpha, \beta}) \cap U \not =\emptyset$, for some clopen set $U$ belonging to the standard base of $\kappa^*$ and $\bigcup_{(\gamma, \delta) <_* (\alpha, \beta)} \mathcal{F}_{\gamma, \delta} \in U$,
			\item  $U_{\alpha, \beta} = f_\alpha^{-1}(U) \subseteq B_{\alpha, \beta}$
			\item  for a decreasing sequence $(V_{\alpha, \beta}^r \colon r \in [\kappa]^{< \omega})$ if each $V^r_{\alpha, \beta} \in \mathcal{F}_{\alpha, \beta}$ then there are $W_{\alpha, \beta}^\nu \in \mathcal{F}_{j(\alpha, \beta)}$, for each $\eta \in \kappa^+$, such that for each $p \in [\kappa^+]^{<\omega}$, $$|\bigcap_{\nu \in p}W_{\alpha, \beta}^\nu\setminus V_{\alpha, \beta}^{\hat{\varphi}(p)}|< \kappa.$$
		\end{itemize}
	\end{itemize} 	
	
	The first and the limit step are done. Assume that we have constructed $\mathcal{F}_{\alpha, \beta}$, $I_{\alpha, \beta}$ and $U_{\alpha, \beta}$. We will show the successor step.  
	
	Assume that  $f_\alpha(B_{\alpha, \beta}) \cap \kappa^* \not = \emptyset$. Otherwise, by $(7)$ we put $$U_{j(\alpha, \beta)} = B_{j(\alpha, \beta)} \cap f^{-1}_\alpha(\{\xi\})$$ for some $\xi \in \kappa$. Then $(ii)$ holds.
	
	Consider two cases.
	
	Case 1. $\alpha+\beta \equiv 0 \ (mod\  2)$. Let  $\mathcal{R} = [\mathcal{F}_{\alpha, \beta}, \{Z_{\alpha, \beta}\}]$ be a filter.
	If $\mathcal{R}$ is proper and 
	$$\{E_t^i \colon t \in [\kappa]^{<\omega}, i \in I_{\alpha, \beta}\} \cup \{A_\mu^i \colon \mu < \kappa^+, i \in I_{\alpha, \beta}\}$$ is independent with respect to $\mathcal{R}, \hat{\varphi}$, then  put 
	then we put
	$$\mathcal{F}_{j(\alpha, \beta)} = \mathcal{F}_{\alpha, \beta}, \ \ I_{j(\alpha, \beta)} = I_{\alpha, \beta}, \ \ U_{j(\alpha, \beta)} = U_{\alpha, \beta}.$$
	
	Otherwise, fix $n \in \omega$, distinct $i_k \in I_{\alpha, \beta}$ and $\hat{\varphi}(p_k)\subseteq t_k,$ for $ k < n$, such that
	$$\kappa \setminus [Z_{\alpha, \beta} \cap \bigcap_{k=0}^{n-1}(\bigcap_{\mu\in p_k}A^{i_k}_{\mu} \cap E^{i_k}_{t_k})]\in \mathcal{F}_{\alpha, \beta}.$$
	Then, put $$\mathcal{F}_{j(\alpha, \beta)} = [\mathcal{F}_{\alpha, \beta}, \{A^{i_k}_{p_k} \colon 0 \leqslant k \leqslant n-1\}, \{E^{i_k}_{t_k} \colon 0 \leqslant k \leqslant n-1\}]$$ $$I_{j(\alpha, \beta)} = I_{\alpha, \beta}\setminus \{i_k \colon 0 \leqslant k \leqslant n-1\}$$
	and
	$$U_{j(\alpha, \beta)} = f^{-1}(U) \subseteq B_{j(\alpha, \beta)},$$
	where $U$ is a clopen set belonging to the standard basis of $\kappa^*$ and $\bigcup_{(\gamma, \delta) <_* j(\alpha, \beta)}\mathcal{F}_{\gamma, \delta} \in U$.
	Then $\kappa\setminus Z_{\alpha, \beta} \in \mathcal{F}_{j(\alpha, \beta)}.$
	
	To show that $(2)$ is fulfilled for $\mathcal{F}_{j(\alpha, \beta)}$ and $I_{j(\alpha, \beta)}$  it is enough to observe that each element of $\mathcal{F}_{j(\alpha, \beta)}$ is of the form 
	$$A \cap (\bigcap_{\mu \in p_k}A^{i_k}_\mu \cap E^{i_k}_{t_k})$$
	for some $A \in \mathcal{F}_{\alpha, \beta} $ and $ 0 \leqslant k \leqslant n-1$. Then,
	$$\bigcap_{\mu \in p_k}A^{i_k}_\mu \cap E^{i_k}_{t_k} \supseteq\bigcap_{k=0}^{n-1}(\bigcap_{\mu\in p_k}A^{i_k}_{\mu} \cap E^{i_k}_{t_k}) $$
	and $$\bigcap_{k=0}^{n-1}(\bigcap_{\mu\in p_k}A^{i_k}_{\mu} \cap E^{i_k}_{t_k}) \in \mathcal{F}_{\alpha, \beta}.$$
	Thus, $(2)$ for $\mathcal{F}_{j(\alpha, \beta)}$ and $I_{j(\alpha, \beta)}$ is fulfilled which follows from $(2)$ for $\mathcal{F}_{\alpha, \beta}$ and $I_{\alpha, \beta}$. 
	
	Case 2. $\alpha +\beta \equiv 1 \ (mod\  2)$.
	
	If $V^r_{\alpha, \beta} \not \in \mathcal{F}_{\alpha, \beta}$ for some $r \in [\kappa]^{<\omega}$ then we put
	$$\mathcal{F}_{j(\alpha, \beta)} = \mathcal{F}_{\alpha, \beta}, \ \ I_{j(\alpha, \beta)} = I_{\alpha, \beta}, \ \ U_{j(\alpha, \beta)} = U_{\alpha, \beta}.$$
	
	Assume that each $V^r_{\alpha, \beta} \in \mathcal{F}_{\alpha, \beta}$. Choose $i \in I_{\alpha, \beta}$. Then for each $\nu < \kappa^+$ consider
	$$W_{\alpha, \beta}^\nu = A^i_\nu \cap \bigcup\{V^t_{\alpha, \beta} \cap E^i_t \colon t \in [\kappa]^{< \omega}\}.$$
	Then for each $p \in [\kappa^+]^{<\omega}$ we have
	$$\bigcap_{\nu\in p} W_{\alpha, \beta}^\nu = \bigcap_{\nu\in p} A^i_\nu \cap \bigcup\{V^t_{\alpha, \beta} \cap E^i_t \colon t \in [\kappa]^{<\omega}\} \subseteq_*$$
	then by $(e)$
	$$\bigcup\{V^t_{\alpha, \beta} \cap E^i_t \colon t \in [\kappa]^{< \omega}\} \cap \{E^i_t \colon t \supseteq \hat{\varphi}(p)\} =$$
	$$\bigcup \{V^t_{\alpha, \beta} \cap E^i_t \colon t \supseteq \hat{\varphi}(p)\} \subseteq $$
	$$\{V^t_{\alpha, \beta} \colon t \supseteq \hat{\varphi}(p)\} = \bigcap_{\nu \in \hat{\varphi}(p)}V^\nu_{\alpha, \beta}.$$
	The last equality is obtained by the monotonicity of $\{V^r_{\alpha, \beta} \colon r \in [\kappa]^{<\omega}\}.$
	
	Now, put 
	$$\mathcal{F}_{j(\alpha, \beta)} = [\mathcal{F}_{\alpha, \beta}, \{W^\nu_{\alpha, \beta} \colon \nu< \kappa^+\}]$$
	$$I_{j(\alpha, \beta)} = I_{\alpha, \beta} \setminus \{i\}$$
	and
	$$U_{j(\alpha, \beta)} = f^{-1}(U) \subseteq B_{j(\alpha, \beta)},$$
	where $U$ is a clopen set belonging to the standard basis of $\kappa^*$ and $\bigcup_{(\gamma, \delta) <_* j(\alpha, \beta)}\mathcal{F}_{\alpha, \beta} \in U$.

	To show that $(2)$ holds for $\mathcal{F}_{j(\alpha, \beta)}$ and $I_{j(\alpha, \beta)}$ it is enough to observe that each element of $\mathcal{F}_{j(\alpha, \beta)}$  is of the form 
	$$A \cap \bigcap_{\nu \in p}W_{\alpha, \beta}^\nu$$
	for some $A\in \mathcal{F}_{\alpha, \beta}$ and $p \in [\kappa^+]^{<\omega}$.
	But 
	$$V^{\hat{\varphi}(p)}_{\alpha, \beta} \cap \bigcap_{\nu \in p} A^i_\nu\cap E^i_{\hat{\varphi}(p)} \subseteq \bigcap_{\nu \in p}W^\nu_{\alpha, \beta}$$
	and since each $V^r_{\alpha, \beta} \in \mathcal{F}_{\alpha, \beta}$ we have that $V^{\hat{\varphi}(p)}_{\alpha, \beta} \in \mathcal{F}_{\alpha, \beta}$.
	Thus, $(2)$ for $\mathcal{F}_{j(\alpha, \beta)}$ and $I_{j(\alpha, \beta)}$ is fulfilled which follows from $(2)$ for $\mathcal{F}_{\alpha, \beta}$ and $I_{\alpha, \beta}$. 
\end{proof}
\\

\begin{theorem}
	Let $\kappa$ be a regular cardinal and let $\hat{\varphi}$  is a $\kappa$-shrinking function.
	Let $X$ be a compact space and let $$\{f_\alpha \colon \alpha< 2^\kappa, f_\alpha \colon X \stackrel{onto}{\longrightarrow} \beta \kappa\}$$ be a family of open mappings. Then there exists a family $\{U_{\alpha, \beta} \colon (\alpha, \beta) \in 2^\kappa\times 2^\kappa\}$ such that for each $x \in \bigcap \{U_{\alpha, \beta} \colon (\alpha, \beta) \in 2^\kappa\times 2^\kappa\}$ and for each $\alpha \in 2^\kappa$  either $f_\alpha(x)$ belongs to $\kappa$ or $f_\alpha(x)$ is a $\kappa$-shrinking point in $\kappa^*$.
\end{theorem}

\begin{proof}
	By Theorem 3, there exists a family $\{U_{\alpha, \beta} \colon (\alpha, \beta) \in 2^\kappa\times 2^\kappa\}$ of non-empty and closed subsets of $X$ and an ultrafilter $\bigcup_{(\alpha, \beta) \in 2^\kappa\times 2^\kappa}\mathcal{F}_{\alpha, \beta}$ of properties $(i)-(iii)$. 
	
	Then for each $x \in \bigcap\{U_{\alpha, \beta} \colon (\alpha, \beta) \in 2^\kappa\times 2^\kappa\}$
and for each $\alpha \in 2^\kappa$, $f_\alpha(x)$ belongs to $\kappa$ whenever
	$f_{\alpha}(B_{\alpha, \beta}) \cap \kappa^* = \emptyset$ or $f_\alpha(x)$    is a $\kappa$-shrinking point in $\kappa^*$
	(of the form $\bigcup_{(\alpha, \beta) \in 2^\kappa \times 2^\kappa}\mathcal{F}_{\alpha, \beta}$) whenever 
	$ f_{\alpha}(B_{\alpha, \beta}) \cap \kappa^* \not = \emptyset$.
\end{proof}
\\

 \section{Remarks and open problems}
 
 Notice that theorems proved in Section 3 are also true if we consider a family of functions $\{f_\alpha \colon \alpha< 2^\lambda, f_\alpha \colon X \stackrel{onto}{\longrightarrow} \beta \kappa\}$, where $\lambda, \kappa$ are regular with $\lambda< \kappa$ and $X$ - a compact space.
 
 We know that in Theorem 2 and in Theorem 4 one cannot substitute "weak $P_{\kappa}$ point" by "weak $P_{\kappa^+}$-point", because in $\beta \kappa$ there are no weak $P_{\kappa^+}$-points, see \cite{BK}. However, the following theorem is an easy consequence of Fact 2, Fact 3 and Theorem 4.
 
 \begin{theorem}
 	Let $\kappa$ be a regular cardinal and let $\hat{\varphi}$  is a $\kappa$-shrinking function.
 	Let $X$ be a compact space and let $$\{f_\alpha \colon \alpha< 2^\kappa, f_\alpha \colon X \stackrel{onto}{\longrightarrow} \beta \kappa\}$$ be a family of open mappings. Then there exists a family $\{U_{\alpha, \beta} \colon (\alpha, \beta) \in 2^\kappa\times 2^\kappa\}$ such that for each $x \in \bigcap \{U_{\alpha, \beta} \colon (\alpha, \beta) \in 2^\kappa\times 2^\kappa\}$ and for each $\alpha \in 2^\kappa$ either $f_\alpha(x)$ belongs to $\kappa$ or $f_\alpha(x)$  is a weak $P_{\kappa}$-point in $\kappa^*$.
 \end{theorem}

We do not know whether theorems proved in Section 3 hold for $\kappa$ singular, although the existence of weak $P_\kappa$-point in $\beta \kappa$ is known also for $\kappa$ singular, see \cite{BK}.

 \noindent
 {\sc Joanna Jureczko}
 \\
 Wroc\l{}aw University of Science and Technology, Wroc\l{}aw, Poland.
 \\
 {\sl e-mail: joanna.jureczko@pwr.edu.pl}
 
\end{document}